\documentclass{aims}

\usepackage{cite}                                       
\usepackage{amsmath, amsbsy, amsfonts, amssymb}
\usepackage{enumerate}
\usepackage{hyperref}

\def\Dscr{\mathcal{D}}
\def\Cscr{\mathcal{C}}

\def\Gscr{\mathcal{G}}

\def\Oscr{\mathcal{O}}
\def\Pscr{\mathcal{P}}

\def\interi{\mathbb{Z}}
\def\reali{{\mathbb{R}}}
\def\complessi{{\mathbb{C}}}
\def\toro{\mathbb{T}}

\let\epsilon\varepsilon
\let\phi\varphi
\def\imunit{{\bf i}}
\def\Lie#1{L_{#1}}

\def\parder#1#2{\frac{\partial#1}{\partial#2}}

\def\typeofarticle{Research Article}
\def\currentvolume{***}
\def\currentissue{***}

\def\ppages{***}
\def\DOI{***}
\def\Received{***}
\def\Accepted{***}
\def\Published{***}

\newtheorem{theorem}{Theorem}[section]
\newtheorem{lemma}{Lemma}[section]

\newtheorem{proposition}{Proposition}[section]

\def\bibref#1{$^{\text{\cite{#1}}}$}

\begin{document}

\title{Kolmogorov variation: KAM with knobs \emph{(\`a la Kolmogorov)}}

\author{Marco Sansottera\affil{1}\corrauth, Veronica Danesi\affil{2}}

\shortauthors{M. Sansottera and V. Danesi}

\address{%
  \addr{\affilnum{1}}{Rocketloop GmbH, Hansaallee 154, 60320 Frankfurt, Germany}
  \addr{\affilnum{2}}{Department of Mathematics, University of Rome “Tor Vergata”, via della Ricerca Scientifica 1, 00133 Rome, Italy.}}
\corraddr{marco@rocketloop.de}

\begin{abstract}
In this paper we reconsider the original Kolmogorov normal form
algorithm~\cite{Kolmogorov-1954} with a variation on the handling of
the frequencies.  At difference with respect to the Kolmogorov
approach, we do not keep the frequencies fixed along the normalization
procedure.  Besides, we select the frequencies of the final invariant
torus and determine \emph{a posteriori} the corresponding starting
ones.  In particular, we replace the classical \emph{translation step}
with a change of the frequencies.  The algorithm is based on the
original scheme of Kolmogorov, thus exploiting the fast convergence of
the Newton-Kantorovich method.
\end{abstract}

\maketitle

\section*{Foreword}
\noindent
The present manuscript marks a first step in answering a question
raised by Prof.~Antonio Giorgilli in 2014 about our recent
result~\cite{GioLocSan-2014} on the construction of lower dimensional
elliptic tori in planetary systems.  The question sounded pretty much
like ``can we fix the final frequencies and determine where we have to
start from?'' Indeed, from the \emph{mathematical} point of view the
result in~\cite{GioLocSan-2014} was satisfactory, we obtained a result
that is valid in measure.  However, Antonio has always pursued
explicit algorithms that can be effectively implemented in order to
study the behavior of a specific dynamical system, like the dynamics
of the solar system or the FPU problem.  Thus, the fact that given a
specific value of the frequencies one does not know if the
corresponding lower dimensional torus exists or not, left us with a
bad taste in our mouth\footnote{As Antonio said, the result was
satisfactory for an \emph{analyst}, not for a \emph{mathematical
physicist}.}.

We therefore believe that the manuscript has an appropriate place in
this Volume in honor of Prof.~Antonio Giorgilli.  We will do our best to
follow the line traced by Antonio and preserve his legacy, always
looking for rigorous \emph{constructive} results\footnote{Where
  \emph{constructive} means that one must be able to \emph{build}:
  (i)~the proof of the theorem, (ii)~the code that implement it;
  (iii)~the computer that run the code; and of course (iv)~the desk
  and the chair where we actually write down the proof.}.

\section{Introduction}
\noindent
The aim of this paper is to reconsider the proof of the Kolmogorov
theorem~\cite{Kolmogorov-1954} with a variation on the handling
of the frequencies.

\subsection{About the genesis of this approach}
The motivation behind the development of this approach has strong
connections with the problem of the persistence of lower dimensional
elliptic invariant tori under sufficiently small perturbations.
Indeed, in \cite{GioLocSan-2014} the authors gave an \emph{almost}
constructive proof of the existence of lower dimensional elliptic tori
for planetary systems, adapting the classical Kolmogorov
normalization algorithm (see also~\cite{SanLocGio-2011}) and a result
of P\"oschel\bibref{Poeschel-1989}, that allows to estimate the
measure of a suitable set of non-resonant frequencies.  The key point
is that both the \emph{internal} frequencies of the torus and the
\emph{transversal} ones vary at each normalization step, and cannot be
kept fixed as in Kolmogorov algorithm.  This makes the accumulation
of small divisors much more tricky to control and, more important, the
result is only valid \emph{in measure} and therefore one cannot know
\emph{a priori} if a specific invariant torus exists or not.

A different approach, based on Lindstedt's series, that allows to
control the frequencies has been proposed
in~\cite{ChrEftBou-2010,ChrEft-2013} in the context of FPU problem.
However, the algorithm has been so far introduced and used, up to our
knowledge, only in a formal way and the literature lacks of rigorous
convergence estimates.  Recently, a comparison of the Lindstedt's method
and the Kolmogorov normal form has been studied in~\cite{EftMas-2020}.

The idea is to overcome the issue of having a result that is valid
only in measure, playing with the frequency like one does with a
control \emph{knob}, hence the title of the paper.  The present work
focuses on full dimensional invariant tori, thus representing a first
step in this direction.  We are well aware that, considering full
dimensional invariant tori, the original Kolmogorov normalization
algorithm allows to have a complete control of the frequencies, which
are kept fixed along the whole normalization procedure.  However,
considering lower dimensional elliptic tori, as explained in detail by
P\"oschel\bibref{Poeschel-1989}, one cannot keep the frequencies
fixed, but have to let them vary.  Thus, as a first result, we decide
to adapt the classical Kolmogorov normalization algorithm in order to
avoid the translation that keeps the frequencies fixed by introducing
a \emph{detuning}\footnote{The detuning can be figured as the action
of turning a control \emph{knob}.}  between the prescribed
\emph{final} frequencies and the corresponding initial ones, to be
determined \emph{a posteriori}. We remark that a similar approach has
been adopted in~\cite{SteLoc-2012}, dealing with an application of the KAM
theorem in dissipative dynamical systems.

Finally, let us stress that our approach (see
also~\cite{Poeschel-2001,Fejoz-2004}), in principle,
also allows to start from a resonant torus that \emph{by construction}
falls into a strongly nonresonant one.

\subsection{KAM theory}

In order to better illustrate our point of view, we briefly recall
here some classical results on KAM theory.  Consider the so-called
fundamental problem of dynamics as stated by Poincar\'e, i.e., a
canonical system of differential equations with Hamiltonian
\begin{equation}
  H(p,q)= H_0(p) + \epsilon H_1(p,q;\epsilon)\ ,
  \label{frm:H}
\end{equation}
where $(p,q)\in\Gscr\times\toro^n$ are action-angle variables, being
$\Gscr\subseteq\reali^n$ an open set and $\epsilon$ is a small
parameter.  The functions $H_0(p)$ and $H_1(p,q;\epsilon)$ are assumed to
be analytic in the variables and in the small parameter, and bounded.
Kolmogorov\bibref{Kolmogorov-1954}, in his seminal paper, that
together with the works of Moser\bibref{Moser-1962} and
Arnold\bibref{Arnold-1963} gave birth to the celebrated KAM theory,
proved the existence of quasi periodic solutions for this Hamiltonian,
with given strongly nonresonant frequencies.

The original idea of Kolmogorov is to select the actions $p^*\in\Gscr$
such that the frequency vector $\omega=\nabla_p{H_0}(p^*)$ satisfies a
Diophantine condition
\begin{equation}
|k\cdot\omega| > \gamma |k|^{-\tau}\ ,\qquad\hbox{for all }k\in\interi^n, k\neq0\ ,
\label{eq:diophantine}
\end{equation}
for some positive $\gamma$ and $\tau\geq n-1$. Hence the term $H_0$
in~\eqref{frm:H} can be expanded in a neighborhood of $p^*$, denoting
again by $p$ the translated actions, and (forgetting the unessential
constant term) the Hamiltonian reads
\begin{equation}
H(p,q)= \omega \cdot p + \Oscr(p^2) + \epsilon H_1(p,q;\epsilon)\ .
  \label{frm:Htrs}
\end{equation}
The Kolmogorov theorem ensures the persistence of the
torus $p=0$ ($p=p^*$ in the original variables) carrying quasi-periodic
solutions with frequencies $\omega$, if $\epsilon$ is small enough and
$H_0(p)$ is nondegenerate.

Let us stress here a technical point.  The role of the nondegeneracy
assumption on $H_0(p)$ is twofold: (i)~it allows to select the desired
frequencies, parameterized by the actions; (ii)~it allows to perform
the \emph{translation step} that keeps the frequency fixed along the
normalization procedure.  However, if the Hamiltonian is already in
the form~\eqref{frm:Htrs} or satisfies the so-called \emph{twistless}
property, i.e., it consists of a sum of a kinetic term, quadratic in
$p$, and of a potential energy, depending only on the angles, it turns
out that the nondegeneracy assumption can be removed, see, e.g.,
\cite{Gallavotti-1994a, GioLoc-1997a}.

Nowadays, the literature about KAM theory is so vast that an
exhaustive list would fill several pages.  Indeed, quoting
P\"oschel\bibref{Poeschel-2001}, \emph{After all, KAM theory is not
  only a collection of specific theorems, but rather a methodology, a
  collection of ideas of how to approach certain problems in
  perturbation theory connected with ``small divisors"}.  Hence, as
this is a paper in honor of Antonio Giorgilli, we have decided to just
mention his main contributions in the field\footnote{For an historical
  account on the role played by Antonio on the development of KAM
  theory in Milan, see~\cite{Galgani-2010} in this same Special
  Issue.}, i.e.,~\cite{BenGalGioStr-1984, Gio-1998, GioLoc-1997a, GioLoc-1997b,
  GioLoc-1999, GioLocSan-2014, GioLocSan-2017, GioLocSan-2009,
  GioMor-1997, LocGio-2000, LocGio-2002, LocGio-2005, MorGio-1995a,
  MorGio-1995b, SanLocGio-2013}.

A final remark is about the so-called \emph{quadratic} (or
superconvergent or Newton-like) method, originally adopted by
Kolmogorov and considered crucial until
Russmann\bibref{Russmann-1989,Russmann-1994} pointed out that a
careful analysis of the accumulation of the small divisors allows to
sharpen some estimates and get rid of it.  Eventually, a proof of
Kolmogorov theorem via classical expansion in a small parameter has
been obtained by Giorgilli and Locatelli in~\cite{GioLoc-1999}.
The approach based on classical expansions allows
to unveil the mechanism of the accumulation of the small divisors and
leads in a natural way to introduce a more relaxed nonresonant
condition for the frequency vector $\omega$, the so-called
$\boldsymbol{\tau}$-condition introduced by Antonio in~\cite{GioMar-2010} and
later adopted in~\cite{GioLocSan-2014, GioLocSan-2015}, precisely
\begin{equation}
  -\sum_{r\geq1}\frac{\log \alpha_r}{r(r+1)} = \Gamma < \infty\ ,
  \qquad\hbox{with}\quad
  \min_{0<|k|\leq rK} |k\cdot\omega| \geq \alpha_r\ ,
  \label{frm:tau}
\end{equation}
where $K$ and $\Gamma$ are two positive constants.  Such a
non-resonance condition is equivalent to the Bruno's one, which is the
weakest one that can be assumed to prove the persistence of invariant
tori (see \cite{Davie-1994,Yoccoz-1988,Yoccoz-1995,GioLoc-1997a}).
Furthermore, the classical approach is the only way to directly
implement KAM theory in practical applications via computer algebra
(see, e.g.,~\cite{GioSan-2011}) and it proved advantageous in
different contexts, e.g., the construction of lower dimensional
elliptic tori in planetary systems
in~\cite{SanLocGio-2011,SanLocGio-2013}, the study of the long term
dynamics of exoplanets in~\cite{LibSan-2013,VolLocSan-2018,
  SanLib-2019}, the investigation of the effective stability in the
spin-orbit problem in~\cite{SanLhoLem-2014,SanLhoLem-2015}, the design
of an a priori control for symplectic maps related to betatronic
motion in~\cite{SanGioCar-2016} and the continuation of periodic
orbits on resonant tori in~\cite{PenSanDan-2018, SanDanPenPal-2020,
  PenDanPal-2021}.

In the present paper, we adopt the original quadratic approach by
Kolmogorov, which turns out to be better suited in order to devise a
normal form algorithm that introduces a detuning of the initial
frequencies that will be determined along the normalization procedure
and complement it with rigorous convergence estimates.

\subsection{Statement of the main result}

Consider a $2n$-dimensional phase space with canonical
action-angle variables $(p,q)\in\Gscr\times\toro^n$, where
$\Gscr\subseteq\reali^n$ is an open set containing the origin.

The Hamiltonian~\eqref{frm:H} is assumed to be a bounded real analytic
function for sufficiently small values of $\epsilon$ and real bounded
holomorphic function of the $(p,q)$ variables in the complex domain
$\Dscr_{\rho_0,\sigma_0} = \Gscr_{\rho_0}\times\toro_{\sigma_0}^{n}$
where $\rho_0$ and $\sigma_0$ are positive parameters,
$\Gscr_{\rho_0}=\bigcup_{p\in\Gscr} \Delta_{\rho_0}(p)$, with
$\Delta_{\rho_0}(p)=\{z\in\complessi^n\colon |p_j-z_j|<\rho_0\}$ and
$\toro_{\sigma_0}^{n}=\{q\in\complessi^n\colon
|\text{Im}(q_j)|<\sigma_0\}$ that are the usual complex extensions of
the real domains.

Given a point $p_0\in\Gscr$, denote by $\omega_0(p_0)\in\reali^n$ the
corresponding frequency vector and expand the Hamiltonian $H_0$ in a
neighborhood of $p_0$, denoting again by $p$ the translated actions
$p-p_0$, precisely
\begin{equation}
  H(p,q)= \omega_0 \cdot p + \Oscr(p^2) + \epsilon H_1(p,q;\epsilon)\ .
  \label{frm:H0trs}
\end{equation}
As remarked in the previous subsection, one can assume a nondegeneracy
condition on $H_0(p)$ so as to ensure that the frequency vector is
parameterized by the actions.  However, if the Hamiltonian is already
in this form, no nondegeneracy assumption is required.

We can now state our main theorem
\begin{theorem}\label{thm:kam-manopole}
  Consider the Hamiltonian~\eqref{frm:H0trs} and pick
  a strongly nonresonant frequency vector $\omega\in\reali^n$
  satisfying the Diophantine condition~\eqref{eq:diophantine} with
  some $\gamma>0$ and $\tau\geq n-1$.  Then there exists a positive
  $\epsilon^*$ such that the following statement holds true: for
  $|\epsilon|<\epsilon^*$ there exist a frequency vector $\omega_0$ and a real analytic near to the
  identity canonical transformation
  $(p,q)=\Cscr^{(\infty)}(p^{(\infty)},q^{(\infty)})$ leading the
  Hamiltonian~\eqref{frm:H0trs} in normal form, i.e.,
  \begin{equation}
    H^{(\infty)} = \omega \cdot p^{(\infty)} + \Oscr({p^{(\infty)}}^2)\ .
    \label{eq:Hkam}
  \end{equation}
\end{theorem}

A more quantitative statement, including a detailed definition of the
threshold on the smallness of the perturbation, is given in
Proposition~\ref{pro:quantitativa}.

A few comments are in order.  At difference with respect to the
original Kolmogorov theorem, we do not keep the frequencies fixed
along the normalization procedure.  The idea, that will be fully
detailed in the next section, is to replace the classical
\emph{translation step} with an \emph{unknown} detuning $\delta\omega$
of the frequencies. Thus, once selected the \emph{final} KAM torus,
the theorem ensures the existence of a \emph{starting} one which is
invariant in the integrable approximation with $\epsilon=0$ and, by
construction, falls into the wanted invariant torus. Let us remark
that in order to apply the Kolmogorov theorem, e.g., for constructing
an invariant torus for a planetary system, it is somehow natural to
determine the final angular velocity vector $\omega$ by using some
numerical techniques like, e.g., Frequency Analysis (see
\cite{Laskar-1994, Laskar-2005}).

\section{Analytic setting and expansion of the Hamiltonian}
We now define the norms we are going to use.  For real vectors
$x\in \reali^n$, we use
$$
|x|=\sum_{j=1}^n |x_j| \ .
$$
For an analytic function $f(p,q)$ with $q\in \toro^n$, we use the
weighted Fourier norm
$$
\|f\|_{\rho,\sigma}=\sum_{k\in \interi^n} |f_k|_{\rho} e^{|k|\sigma} \ ,
$$
with 
$$
|f_k|_{\rho}=\sup_{p} |f_k(p)| \ .
$$
We introduce the classes of
functions $\Pscr_{l}$, with integers $l\geq0$, such that $g\in\Pscr_{l}$ can be
written as
$$
g(p,q) = \sum_{|m|=l} \sum_{k} c_{m,k} p^m e^{\imunit k \cdot q}\ ,
$$ with $c_{m,k}\in\complessi$. For consistency reasons, we also set
$\Pscr_{-1}=\{0\}$.  Finally, we will also omit the
dependence of the functions from the variables, unless it has some
special meaning.

The Hamiltonian~\eqref{frm:H0trs}, expanded in power series of the actions $p$, reads
\begin{equation}
  H(p,q) = \omega_0 \cdot p + \sum_{l\geq0} h_{l}
  \label{frm:Hexp0}
\end{equation}
where $h_{l}\in\Pscr_{l}$ are bounded as
\begin{equation}
\|h_0\|_{\rho,\sigma} \leq \epsilon E\ ,\quad
\|h_1\|_{\rho,\sigma} \leq \frac{\epsilon E}{2}
\qquad\hbox{and}\qquad
\|h_{l}\|_{\rho,\sigma} \leq \frac{E}{2^l}\quad\hbox{for }l\geq2\ .
\label{eq:stime-h_l}
\end{equation}
provided $\rho_0\leq 1/4$, with $E=2^{n-1} E_0$ where
\begin{equation}
E_0 = \max \left( \sup_{p\in\Delta_{\rho_0}}|H_0(p)|\ ,\ \sup_{(p,q)\in\Dscr_{\rho_0,\sigma_0}}|H_1(p,q;\epsilon)|\right)\ .
\label{frm:E0-G0}
\end{equation}

\section{Formal algorithm}\label{sec:formale}
\noindent
We present in this section the algorithm leading the
Hamiltonian~\eqref{frm:Hexp0} in normal form.  The procedure is
described here from a purely formal point of view, while the study of
the convergence is postponed to the next section.

First we introduce the unknown detuning $\delta\omega$ and rewrite the
Hamiltonian as
\begin{equation}
  H(p,q)= \omega\cdot p + \delta\omega\cdot p + \sum_{l\geq0} h_{l}(p,q)\ ,
  \label{frm:H0}
\end{equation}
with $h_{l}\in\Pscr_{l}$.  Let us stress again that the quantity
$\delta\omega$ is unknown and will be determined at the end of the
normalization procedure.

As in the original Kolmogorov proof scheme, the algorithm consists in
iterating infinitely many times a single normalization step: starting
from $H$, we apply two near to the identity canonical transformations
with generating functions $\chi_{0}(q)$ and $\chi_{1}(p,q)$, i.e.,
$$
H' = \exp(L_{\chi_1})\circ\exp(L_{\chi_0}) H\ .
$$
The generating functions
are determined in order to kill the unwanted terms $h_{0}(q)$ and
$h_{1}(p,q)$.  At difference with respect to the
original approach designed by Kolmogorov we do not introduce a
translation of the actions $p$, since we do not keep fixed the
\emph{initial} frequency $\omega_0$. Indeed, in our algorithm
the role of the translation step is played by the detuning of the
frequency $\delta\omega$.

The functions $\chi_{0}(q)$ and $\chi_{1}(p,q)$ are
determined by solving
\begin{align}
  &L_{\chi_0} \omega \cdot p + h_0 = 0\ ,
  \label{frm:eqom1}\\
  &L_{\chi_1} \omega \cdot p + \sum_{s\geq0}\frac{1}{s!}\,\Lie{\chi_0}^s h_{s+1} =
  \sum_{s\geq0}\frac{1}{s!}\,\langle \Lie{\chi_0}^s h_{s+1} \rangle_{q} \ ,
  \label{frm:eqom2}
\end{align}
where $\langle \cdot \rangle_{q}$ denotes the average with respect to
the angles $q$.

First, considering the Fourier expansion of $h_0$, and neglecting the constant term, one has
$$
h_0(q) = \sum_{k\neq0} c_{k} e^{\imunit  k\cdot q}\ ,
$$
and can easily check
that the solution of~\eqref{frm:eqom1} is given by
$$
\chi_0(q) =  \sum_{k\neq0} \frac{c_{k}}{\imunit k\cdot\omega} e^{\imunit k \cdot q}\ .
$$
The intermediate Hamiltonian $\hat{H}=\exp(L_{\chi_0})H$ reads
\begin{equation}
  \hat{H}(p,q)= \omega\cdot p + \delta\omega'\cdot p
  + \sum_{l\geq0} \hat{h}_{l}(p,q) \ ,
  \label{frm:hat{H}^{(r)}}
\end{equation}
with
\begin{equation}
  \begin{aligned}
    \delta\omega'\cdot{p}&=\delta\omega \cdot{p}+
  \sum_{s=0}^{\infty}\frac{1}{s!}\,\langle{\Lie{\chi_0}^s h_{s+1}}\rangle_q\ ,\\
    \hat h_0 &=\Lie{\chi_0}\Big(\delta\omega'\cdot p -
 \sum_{s=1}^{\infty}\frac{1}{s!}\,\langle\Lie{\chi_0}^s h_{s+1}\rangle_q\Big)\\
 & \quad
 + \Lie{\chi_0} \Big( h_{1} - \langle h_{1}\rangle_q \Big)
 +\sum_{s=2}^{\infty}\frac{1}{s!}\,\Lie{\chi_0}^s h_s\ ,\\
 \hat{h}_1 &=\sum_{s=0}^{\infty}\frac{1}{s!}\,\Lie{\chi_0}^s h_{s+1}
 +\left(\delta\omega-\delta\omega'\right)\cdot{p}\ ,\\
 \hat{h}_l &=
\sum_{s=0}^{\infty}\frac{1}{s!}\,\Lie{\chi_0}^s h_{s+l}\ ,
&\quad\hbox{for }l\ge 2\ .
  \end{aligned}
  \label{frm:hat{H}^{(r)}-functions}
\end{equation}
where the unessential constant term
$\langle h_{0} \rangle_q$ has been neglected in the expression
above.

Second, considering the Fourier expansion
$$
\hat{h}_1(p,q) = \sum_{k\neq 0} \hat{c}_{k}(p) e^{\imunit k \cdot q}\ ,
$$
one can easily check that the solution of~\eqref{frm:eqom2} is given by
$$
\chi_1(p,q) =  \sum_{k\neq 0} \frac{\hat{c}_{k}(p)}{\imunit k \cdot \omega} e^{\imunit k \cdot q}
$$

We complete the normalization step by computing the Hamiltonian
$H'=\exp(L_{\chi_1})\hat{H}$ that takes the
form~\eqref{frm:H0} with $\delta\omega'$ as in~\eqref{frm:hat{H}^{(r)}-functions} and
\begin{equation}
  \begin{aligned}
    h'_0 &=\sum_{s=0}^{\infty}\frac{1}{s!}\,\Lie{\chi_1}^s \hat{h}_0\ , \\
    h'_1 &=\sum_{s=1}^{\infty}\frac{s}{(s+1)!}\,\Lie{\chi_1}^s \hat{h}_{1}
    +\sum_{s=1}^{\infty}\frac{1}{s!}\,\Lie{\chi_1}^s{\delta\omega'}\cdot p\ ,&\\
    h'_l &=\sum_{s=0}^{\infty}\frac{1}{s!} \Lie{\chi_1}^s \hat{h}_{l}
&\quad\hbox{for } l\ge 2\ .
  \end{aligned}
  \label{frm:H^{(r)}-functions}
\end{equation}
The justification of the
formul{\ae}~\eqref{frm:hat{H}^{(r)}-functions}
and~\eqref{frm:H^{(r)}-functions} is just a matter of straightforward
computations, exploiting~\eqref{frm:eqom1} and~\eqref{frm:eqom2}.

\section{Quantitative estimates}\label{sec:stime}
\noindent
In this section, we translate our formal algorithm into a recursive
scheme of estimates on the norms of the functions.  This essentially
requires to bound the norm of the Lie series.  In order to shorten the
notation, we will replace $|\cdot|_{\alpha(\rho,\sigma)}$ by
$|\cdot|_\alpha$ and $\|\cdot\|_{\alpha(\rho,\sigma)}$ by
$\|\cdot\|_\alpha$, being $\alpha$ any positive real number.  The
useful estimates are collected in the following statements.

\begin{lemma}\label{lem:stima-derivata-Lie}
Let $f$ and $g$ be analytic respectively in
$\Dscr_{1}$ and $\Dscr_{(1-d')}$ for some $0\le
d'<1$ with finite norms $\|f\|_{1}$ and
$\|g\|_{1-d'}$. Therefore,
\begin{enumerate}[i.]
\item for $0< d< 1$ and for $1\le j\le n$ we have
\begin{equation}
  \left\|\parder{f}{p_j}\right\|_{(1-d)} \le
  \frac{1}{d\rho}\|f\|_{1}\ ,\quad
  \left\|\parder{f}{q_j}\right\|_{(1-d)} \le
  \frac{1}{ed\sigma}\|f\|_{1}\ ;                                        
  \label{frm:stima-der}
\end{equation}
\item for $0< d< 1-d'$ we have
\begin{equation}
  \|\{f,g\}\|_{(1-d'-d)} \le \frac{2}{ed(d+d')\rho\sigma} \|f\|_{1}\|g\|_{(1-d')} \ .
  \label{frm:poisson}                                                           
  \end{equation}
\end{enumerate}
\end{lemma}

\begin{lemma}\label{lem:stima-termini-serie-Lie}
  Let $d$ and $d^{\prime}$ be real numbers such that $d>0\,$,
  $d^{\prime}\ge 0$ and $d+d^{\prime}<1\,$; let $X$ and $g$ be two
  analytic functions on $\Dscr_{(1-d^{\prime})}$ having
  finite norms $\|X\|_{1-d^{\prime}}$ and $\|g\|_{1-d^{\prime}}\,$,
  respectively. Then, for $j\ge 1$, we have
\begin{equation}
\frac{1}{j!}\left\|L^{j}_{X}g\right\|_{1-d-d^{\prime}}
\le\frac{1}{e^{2}}
\left(\frac{2e}{\rho\sigma}\right)^{j}
\frac{1}{d^{2j}}\|X\|^{j}_{1-d^{\prime}}\|g\|_{1-d^{\prime}}\ .
\label{frm:stimalie}
\end{equation}
\end{lemma}
The proofs of these lemmas are straightforward
and can be found, e.g., in~\cite{Gio-2003}.

We are now ready to write the statement of Theorem 8.1 in a more detailed form.

\begin{proposition}\label{pro:quantitativa}
  Consider the Hamiltonian~\eqref{frm:H0} and assume the following hypotheses:
  \begin{enumerate}[(i)]
  \item $h_l$, for $l\geq0$, satisfy~\eqref{eq:stime-h_l};
  \item $\omega\in\reali^n$ satisfy the
    Diophantine condition~\eqref{eq:diophantine} with some $\gamma>0$
    and $\tau\geq n-1$.
  \end{enumerate}
 Then, there exists a positive $\epsilon^*$ depending on $n$, $\tau$,
 $\gamma$, $\rho$ and $\sigma$ such that for $|\epsilon|<\epsilon^*$ and $\delta\leq1/8$ there
 exists a real analytic near to the identity canonical transformation
 $(p,q)=\Cscr(p',q')$ satisfying
 \begin{equation}
   |p_j-p'_j| \leq \delta^{\tau+3}\rho\ ,\quad
   |q_j-q'_j| \leq \delta^{\tau+3}\sigma\ ,\qquad
   j=1\,,\ \ldots\,,\ n\ ,
 \end{equation}
for all $(p',q')\in\Dscr_{1-4\delta}$ which gives the Hamiltonian the
Kolmogorov normal form~\eqref{eq:Hkam}. Moreover, the detuning is bounded as
$$
\|\delta\omega \cdot p\|_{\frac{1}{2}} \leq \frac{E}{2} \delta^{2\tau+4}\ .
$$
\end{proposition}

The proof of this Proposition is given in the next two
subsections. Indeed, it is divided in two parts: first the quantitative
analytic estimates for a single step are obtained in the so-called
Iterative Lemma, and finally the convergence of the infinite sequence
of iterations is proved.

\subsection{The Iterative Lemma}
The aim of this subsection is to translate the algorithm of
Section~\ref{sec:formale} into a scheme of estimates for the norms of
all functions involved.

\begin{lemma}\label{lem:iterativo}
Let $H$ be as in~\eqref{frm:H0} and assume that the hypotheses
(i)--(ii) of Proposition~\ref{pro:quantitativa} hold true.  Let
$\delta\leq1/8$ and $\rho^*$, $\sigma^*$ be positive constants
satisfying
$$
(1-4\delta) \rho \geq\rho^*
\quad\hbox{and}\quad
(1-4\delta) \sigma \geq\sigma^*\ .
$$
Then there exists a positive constant $\Lambda = \Lambda(n, \tau,
\gamma, \rho^*, \sigma^*)$ such that the following holds true: if
\begin{equation}
\frac{\Lambda}{\delta^{3\tau+6}}\epsilon \leq 1\ ,
\end{equation}
assuming that the following ``a priori'' bound on the detunings holds true,
$$
\|\delta\omega'\cdot p\|_{1-\delta} \leq \frac{\epsilon}{2 \delta^{\tau+2}}\ ,
$$
then there exists a canonical transformation $(p,q)=\Cscr(p',q')$ satisfying
\begin{equation}
  \begin{aligned}
   |p_j-p'_j| &\leq \frac{\Lambda\epsilon}{\delta^{3\tau+6}} \delta^{\tau+3} \rho\leq \delta^{\tau+3}\rho\ ,\\
   |q_j-q'_j| &\leq \frac{\Lambda\epsilon}{\delta^{3\tau+6}} \delta^{\tau+3} \sigma\leq \delta^{\tau+3}\sigma\ ,\qquad
   j=1\,,\ \ldots\,,\ n\ ,
   \label{eq:stima_deformazione}
  \end{aligned}
\end{equation}
for all $(p',q')\in\Dscr_{1-4\delta}$, which brings the Hamiltonian in
the Kolmogorov normal form~\eqref{eq:Hkam} with the same $\omega$ and
with new functions $\delta\omega'\cdot p$ and $h'_{l}$, for $l\geq1$,
satisfying the hypotheses (i)--(ii) with new positive constants
$\epsilon'$, $\rho'$, $\sigma'$ given by
$$
\epsilon' = \frac{\Lambda}{\delta^{3\tau+6}}\epsilon^2\ ,\quad
\rho' = (1-4\delta) \rho\quad\hbox{and}\quad
\sigma' = (1-4\delta) \sigma\ .
$$
Furthermore, the variation of the detuning frequency vector is bounded as follows
$$
\|\delta\omega' \cdot p - \delta\omega\cdot p\|_{1-2\delta} \leq
\frac{\Lambda\epsilon}{\delta^{3\tau+6}} \delta^{2\tau+4} \frac{E}{2}\leq
\delta^{2\tau+4} \frac{E}{2}\ .
$$
\end{lemma}

A crucial role in the proof of the Iterative Lemma is played by the
control of the accumulation of the small divisors.  This topic has
been deeply investigated by Antonio Giorgilli, see, e.g., \cite{Gio-2022}, Section 8.2.4.

We now collect all the estimates that allow to prove
Lemma~\ref{lem:iterativo}.  Recalling the Diophantine
condition~\eqref{eq:diophantine}, the elementary inequality $|k|^\tau
e^{-|k|\delta\sigma}\le\left(\frac{\tau}{e\delta\sigma}\right)^{\tau}$
allows us to easily bound the generating function $\chi_{0}$ as
$$
\|\chi_{0}\|_{1-\delta} \leq \frac{1}{\gamma}
\left( \frac{\tau}{e \delta \sigma} \right)^\tau \epsilon E
\leq \frac{K_1}{\delta^\tau}\epsilon\ ,\qquad
K_1=\frac{1}{\gamma}\left( \frac{\tau}{e \sigma} \right)^\tau E\ .
$$

It is now convenient to provide some useful estimates to bound the
terms appearing in $\hat{H}$. Assuming the smallness condition on $\epsilon$
$$
\frac{2e K_1 \epsilon}{\delta^{\tau+2}\rho\sigma} \leq \frac{1}{2}\ ,
$$
we easily get
\begin{equation}
  \begin{aligned}
  \| \Lie{\chi_0}h_1 \|_{1-2\delta} &\leq
  \frac{1}{e^2} \left( \frac{2e}{\delta^2\rho\sigma}
  \frac{K_1}{\delta^\tau}\epsilon \right) \frac{\epsilon E}{2}
  \leq \frac{K_2}{\delta^{\tau+2}} \frac{1}{2}\epsilon^2\ ,&\quad 
  K_2=\frac{2 K_1 E}{e\rho\sigma}\ ,\\
  \| \Lie{\chi_0}h_{l+1} \|_{1-2\delta} &\leq
  \frac{1}{e^2} \left( \frac{2e}{\delta^2\rho\sigma}
  \frac{K_1}{\delta^\tau}\epsilon \right) \frac{E}{2^{l+1}}
  \leq \frac{K_2}{\delta^{\tau+2}} \frac{1}{2^{l+1}}\epsilon\ ,\\
  \sum_{s\geq2} \frac{1}{s!} \| \Lie{\chi_0}^s h_{l+s} \|_{1-2\delta} &\leq
  \sum_{s\geq2} \frac{1}{e^2} \left( \frac{2e}{\delta^2\rho\sigma}
  \frac{K_1}{\delta^\tau}\epsilon \right)^s \frac{E}{2^{l+s}}\\
  &\leq\frac{1}{e^2} \left( \frac{2e}{\delta^2\rho\sigma}
  \frac{K_1}{\delta^\tau}\epsilon \right)^2\frac{E}{2^{l+2}}
  \sum_{s\geq0} \left( \frac{e K_1 \epsilon}{\delta^{\tau+2}\rho\sigma}\right)^s&\\
  &\leq \frac{K_3}{\delta^{2\tau+4}} \frac{1}{2^{l+2}}\epsilon^2\ ,
  &\quad K_3=\frac{2^3 K_1^2 E}{\rho^2\sigma^2}\ ,\\
  \sum_{s\geq2} \frac{s}{s!} \| \Lie{\chi_0}^s h_{s} \|_{1-2\delta} &\leq
  \sum_{s\geq2} \frac{s}{e^2} \left( \frac{2e}{\delta^2\rho\sigma}
  \frac{K_1}{\delta^\tau}\epsilon \right)^s \frac{E}{2^{s}}\\
  &\leq\frac{1}{e^2} \left( \frac{2e}{\delta^2\rho\sigma}
  \frac{K_1}{\delta^\tau}\epsilon \right)^2\frac{E}{2^{2}}
  \sum_{s\geq0} (s+2)\left( \frac{e K_1 \epsilon}{\delta^{\tau+2}\rho\sigma}\right)^s&\\
  &\leq \frac{K_4}{\delta^{2\tau+4}} \frac{1}{2^{2}}\epsilon^2\ ,
  &\quad K_4=\frac{24 K_1^2 E}{\rho^2\sigma^2}\ ,
\end{aligned}
\end{equation}
where, in the last two inequalities, we used the well known sums
$$
\sum_{s\geq0} x^s =\frac{1}{1-x} \leq 2 \quad\hbox{and}\quad
\sum_{s\geq1} s x^s =\frac{x}{(1-x)^2} \leq 2\ , \quad\hbox{for }|x|\leq\frac{1}{2}\ .
$$

We now estimate the difference between the detunings, precisely,
\begin{equation}
  \begin{aligned}
    \|\delta\omega' \cdot p - \delta\omega\cdot p\|_{1-2\delta} &\leq
    \sum_{s\geq0} \frac{1}{s!} \| \Lie{\chi_0}^s h_{l+1} \|_{1-2\delta} \\
    &\leq\frac{\epsilon E}{2} +
    \frac{K_2}{\delta^{\tau+2}}\frac{1}{2^2}\epsilon+
    \frac{2^3 K_1^2 E}{\delta^{2\tau+4}\rho^2\sigma^2}\frac{1}{2^3}\epsilon^2\\
    &\leq\frac{\epsilon E}{2} +
    \frac{K_2}{\delta^{\tau+2}}\frac{1}{2^2}\epsilon+
    \left( \frac{2eK_1\epsilon}{\delta^{\tau+2}\rho\sigma} \right)
    \frac{K_2}{\delta^{\tau+2}}\frac{1}{2^2}\epsilon&\\
    &\leq\frac{K_5}{\delta^{\tau+2}}\frac{\epsilon }{2}\ , &\quad K_5=E+K_2\ .
  \end{aligned}
  \label{frm:stima-detunings}
\end{equation}

We now bound the term appearing in~\eqref{frm:hat{H}^{(r)}}.  The norm
of the function $\hat h_0$ is bounded as
$$
\begin{aligned}
\|\hat h_0\|_{1-2\delta} &\leq
\biggl\| \Lie{\chi_0}\Bigl(\delta\omega'\cdot p -
\sum_{s=1}^{\infty}\frac{1}{s!}\,\langle\Lie{\chi_0}^s h_{s+1}\rangle_q\Bigr)\\
&\qquad
+\Lie{\chi_0} \Big( h_{1} - \langle h_{1}\rangle_q \Big)
+\sum_{s=2}^{\infty}\frac{1}{s!}\,\Lie{\chi_0}^s h_s \biggr\|_{1-2\delta}\\
&\leq \bigl\| \Lie{\chi_0}\delta\omega'\cdot p\bigr\|_{1-2\delta}
+ \biggl\| \sum_{s\geq2}\frac{s}{s!}\Lie{\chi_0}^s h_{s}\biggr\|_{1-2\delta}\\
&\qquad
+ \|\Lie{\chi_0} h_{1} \|_{1-2\delta}+ \biggl\| \sum_{s\geq2}\frac{1}{s!}\Lie{\chi_0}^s h_{s}\biggr\|_{1-2\delta}\\
&\leq\frac{1}{e^2}\left( \frac{2e}{\delta^2\rho\sigma} \frac{K_1}{\delta^\tau}\epsilon\right)
\frac{\epsilon}{2 \delta^{\tau+2}}
+\frac{K_4}{\delta^{2\tau+4}}\frac{1}{2^2}\epsilon^2\\
&\qquad
+\frac{K_2}{\delta^{\tau+2}}\frac{1}{2}\epsilon^2
+\frac{K_3}{\delta^{2\tau+4}}\frac{1}{2^2}\epsilon^2\\
&\leq\frac{K_6}{\delta^{2\tau+4}}\epsilon^2\ ,
&\quad K_6=\frac{K_2}{2E}+\frac{K_4}{2^2}+\frac{K_2}{2}+\frac{K_3}{2^2}\ ,\\
\end{aligned}
$$
while the norm of $\hat h_1$ satisfies
$$
\begin{aligned}
\|\hat h_1\|_{1-2\delta} &\leq
\sum_{s\geq0}\frac{1}{s!} \| \Lie{\chi_0}^s h_{s+1}\|_{1-2\delta}\\
&\leq\frac{\epsilon E}{2}
+\frac{K_2}{\delta^{\tau+2}}\frac{1}{2^2}\epsilon
+\frac{K_3}{\delta^{2\tau+4}}\frac{1}{2^3}\epsilon^2\\
&\leq\frac{\epsilon E}{2}
+\frac{K_2}{\delta^{\tau+2}}\frac{1}{2^2}\epsilon
+\left( \frac{2eK_1\epsilon}{\delta^{\tau+2}\rho\sigma} \right)
    \frac{K_2}{\delta^{\tau+2}}\frac{1}{2^2}\epsilon&\\
    &\leq\frac{K_5}{\delta^{\tau+2}}\frac{\epsilon}{2}\ .\\
\end{aligned}
$$
Finally, $\hat h_l$, for $l\geq2$, one has
$$
\begin{aligned}
\|\hat h_l\|_{1-2\delta} &\leq
\sum_{s\geq0}\frac{1}{s!} \| \Lie{\chi_0}^s h_{s+l}\|_{1-2\delta}\\
&\leq
\frac{ E}{2^l}
+\frac{K_2}{\delta^{\tau+2}}\frac{1}{2^{l+1}}\epsilon
+\frac{K_3}{\delta^{2\tau+4}}\frac{1}{2^{l+2}}\epsilon^2\\
&\leq\frac{ E}{2^l}
+\left( \frac{2eK_1\epsilon}{\delta^{\tau+2}\rho\sigma} \right)
    \frac{E}{2 e^2}\frac{1}{2^l}
+\left( \frac{2eK_1\epsilon}{\delta^{\tau+2}\rho\sigma} \right)^2
    \frac{E}{2 e^2}\frac{1}{2^l}\\
&\leq\frac{K_7}{2^l}\ ,& K_7=E+\frac{E}{e^2}\ .
\end{aligned}
$$

This concludes the estimates for the first half of the normalization
step.

Exploiting again the Diophantine condition~\eqref{eq:diophantine} we
easily bound the generating function $\chi_{1}$ as
$$
\|\chi_{1}\|_{1-3\delta} \leq \frac{1}{\gamma}
\left( \frac{\tau}{e \delta \sigma} \right)^\tau
\frac{K_5}{\delta^{\tau+2}}\frac{1}{2}\epsilon
\leq \frac{K_8}{\delta^{2\tau+2}}\frac{1}{2}\epsilon\ ,\qquad
K_8=\frac{1}{\gamma}\left( \frac{\tau}{e \sigma} \right)^\tau K_5\ .
$$

Assuming the smallness condition
$$
\frac{2e \|\chi_1\|_{1-3\delta}}{\delta^2 \rho \sigma} \leq \frac{1}{2}\ ,
$$
that can be written as
$$
\frac{2e K_8\epsilon}{\delta^{2\tau+4} \rho \sigma} \leq 1\ ,
$$
we now bound the terms appearing in~\eqref{frm:H^{(r)}-functions}.  The norm
of the function $h'_0$ is bounded as
$$
\begin{aligned}
\|h'_0\|_{1-4\delta} &\leq
\sum_{s\geq0}\frac{1}{e^2} \left( \frac{2e}{\delta^2\rho\sigma}
\frac{K_8}{\delta^{2\tau+2}} \frac{1}{2} \epsilon\right)^s
\frac{K_6}{\delta^{2\tau+4}}\epsilon^2\\
&\leq
\frac{K_6}{e^2 \delta^{2\tau+4}}\epsilon^2
\sum_{s\geq0} \left( 
\frac{e K_8 \epsilon}{\delta^{2\tau+4}\rho\sigma}\right)^s\\
&\leq
\frac{K_{9}}{\delta^{2\tau+4}}\epsilon^2\ ,
&\quad K_{9} = \frac{2 K_6}{e^2}\ .
\end{aligned}
$$

Similarly we get
$$
\begin{aligned}
\|h'_1\|_{1-4\delta} &\leq
\sum_{s\geq1}\frac{1}{s!} \left\|
\Lie{\chi_1}^s \hat{h}_1 \right\|_{1-4\delta}+
\sum_{s\geq1}\frac{1}{s!} \left\|
\Lie{\chi_1}^s \delta\omega' \cdot p \right\|_{1-4\delta}\\
&\leq
\sum_{s\geq1}\frac{1}{e^2}
\left( \frac{2e}{\delta^2\rho\sigma}
\frac{K_8}{\delta^{2\tau+2}} \frac{1}{2} \epsilon\right)^s
\left( \frac{K_5}{\delta^{\tau+2}} \frac{1}{2} \epsilon+
\frac{1}{\delta^{\tau+2}} \frac{1}{2} \epsilon
\right)\\
&\leq
\frac{1}{e^2}
\left( \frac{2e}{\delta^2\rho\sigma}
\frac{K_8}{\delta^{2\tau+2}} \frac{1}{2} \epsilon\right)
\left( \frac{K_5}{\delta^{\tau+2}} \frac{1}{2} \epsilon+
\frac{1}{\delta^{\tau+2}} \frac{1}{2} \epsilon
\right)\\
&\qquad\sum_{s\geq0}
\left( \frac{2e}{\delta^2\rho\sigma}
\frac{K_8}{\delta^{2\tau+2}} \frac{1}{2} \epsilon\right)^s\\
&\leq
\frac{K_8}{e \delta^{2\tau+4} \rho \sigma} \left( \frac{K_5+1}{\delta^{\tau+2}} \right)2 \frac{\epsilon^2}{2}\\
&\leq \frac{K_{10}}{\delta^{3\tau+6}} \frac{\epsilon^2}{2}\ ,
& K_{10}=\frac{2(K_5+1)K_8}{e\rho\sigma}\ .
\end{aligned}
$$

For $l\geq2$, we bound the norm of $h'_l$ as
$$
\begin{aligned}
  \|h'_l\|_{1-4\delta} &\leq
  \sum_{s\geq0} \frac{1}{s!} \left\|\Lie{\chi_1}^s \hat{h}_l \right\|_{1-4\delta}\\
  &\leq \sum_{s\geq0} \frac{1}{e^2}
  \left( \frac{2e}{\delta^2\rho\sigma}
  \frac{K_8}{\delta^{2\tau+2}} \frac{1}{2} \epsilon\right)^s \frac{K_7}{2^l}\\
  &\leq \frac{K_{11}}{2^l}\ ,&\quad K_{11}=\frac{2K_7}{e^2}\ .
\end{aligned}
$$

To finish, we need to provide the convergence of the near to the
identity change of coordinates.

The first change of coordinates is bounded as follows
$$
\begin{aligned}
  \exp(\Lie{\chi_0})\hat{p}_j &= \hat{p}_j+\parder{\chi_0}{q_j}\bigg|_{(\hat{q},\hat{p})}\ ,\\
  \exp(\Lie{\chi_0})\hat{q}_j &= \hat{q}_j\ ,
\end{aligned}
$$
The second change of coordinates is bounded as 
$$
\begin{aligned}
  \|\exp(\Lie{\chi_1})p'-p'\|_{1-4\delta} &\leq
  \sum_{s\geq1}\frac{1}{s!} \|\Lie{\chi_1}^sp' \|_{1-4\delta}\\
  &\leq
  \sum_{s\geq1}\frac{1}{s!} \frac{(s-1)!}{e^2}
  \left( \frac{2e}{\delta^2\rho\sigma} \|\chi_1\|_{1-4\delta} \right)^{s-1} \|\Lie{\chi_1}p' \|_{1-4\delta}\\
  &\leq \frac{\|\chi_1\|}{e^3\delta\sigma}
  \sum_{s\geq1}\frac{1}{s}\left( \frac{2e}{\delta^2\rho\sigma} \|\chi_1\|_{1-4\delta} \right)^{s-1}\\
  &\leq \frac{\delta\rho}{2}
  \left(\frac{1}{e^4}\frac{2eK_8\epsilon}{\delta^{2\tau+4}\rho\sigma}\right)\ ,
\end{aligned}
$$
and similar computations give
\begin{equation}\label{stima-q}
\|\exp(\Lie{\chi_1})q'-q'\|_{1-4\delta} \leq \frac{\delta\sigma}{2}
\left(\frac{1}{e^3}\frac{2eK_8\epsilon}{\delta^{2\tau+4}\rho\sigma}\right)\ .
\end{equation}
Combining these bounds we eventually get
\begin{equation}\label{stima-p}
\begin{aligned}
\|p'-p\|_{1-4\delta} &\leq
\frac{\delta\rho}{2}\frac{1}{e^4}\frac{2eK_8\epsilon}{\delta^{2\tau+4}\rho\sigma}
+\frac{1}{e\delta\sigma}\frac{K_1}{\delta^\tau}\epsilon\\
&\leq \frac{\delta\rho}{2}
\left( \frac{1}{e^4}\frac{2eK_8}{\delta^{2\tau+4}\rho\sigma}
+\frac{2K_1}{e\delta^{\tau+2}\rho\sigma}\right)\epsilon\ .
\end{aligned}
\end{equation}

In order to conclude the proof, we now collect all the estimates.  We define $\Lambda$ as
$$
\Lambda=\max\left( 1,\ K_j \hbox{ for }j=1,\ldots,11\ ,\ \frac{2eK_1}{\rho\sigma}\ ,\ \frac{2eK_8}{\rho\sigma}  \right)\ .
$$
Let us stress that $\Lambda$ depends only on $\tau$, $\gamma$,
$\rho^*$, $\sigma^*$ and $n$ (implicitly via $\tau$).  Thus all the
convergence conditions are summarized by
$$
\frac{\Lambda}{\delta^{3\tau+6}}\epsilon \leq 1\ ,
$$
and trivial computations conclude the proof of Lemma~\ref{lem:iterativo}.

\subsection{Conclusion of the Proof}
By repeated application of the Iterative Lemma, we construct an
infinite sequence $\{\hat{\Cscr}^{(k)}\}_{k\geq1}$ of near the identity
canonical transformations
$$
(p^{(k-1)},q^{(k-1)})=\hat\Cscr^{(k)}(p^{(k)},q^{(k)})\ ,
$$
where the upper index labels the coordinates at the $k$-th step.  This
introduces a sequence $\{H^{(k)}\}_{k\geq1}$ of Hamiltonians, where
$H^{(0)}=H$ is the original one, satisfying
\begin{align}
  \epsilon_k &= \frac{\Lambda}{\delta_k^{3\tau+6}}\epsilon_{k-1}^2\ ,\label{eq:epsilon_k}\\
  \rho_k &= (1-4\delta_k) \rho_{k-1}\ ,\label{eq:rho_k}\\
  \sigma_k &= (1-4\delta_k) \sigma_{k-1}\ .\label{eq:sigma_k}
\end{align}  
These sequences depend on the arbitrary sequence
$\{\delta_k\}_{k\geq1}$, that must be chosen so that for every $k$ one
has $\delta_k\leq1/8$ and
\begin{align}
  \frac{\Lambda}{\delta_k^{3\tau+6}}\epsilon_{k-1}&\leq1\ ,\label{eq:eps^2}\\
  (1-4\delta_k) \rho_{k-1}&\geq\rho^*>0\ ,\label{eq:rho*}\\
  (1-4\delta_k) \sigma_{k-1}&\geq\sigma^*>0\ .\label{eq:sigma*}
\end{align}  

Let us now make a choice of the parameters\footnote{The choice is
rather arbitrary, see~\cite{Gio-2022}, footnote 6, chapter 8.},
precisely
$$
\epsilon_k = \epsilon_0^{k+1}\qquad\hbox{and}\qquad \delta_k=\frac{1}{\alpha^k}\ ,
$$
where $\alpha$ is real positive constant to be determined.

Let's start with~\eqref{eq:eps^2}, that reads
\begin{equation}
\Lambda (\alpha^{3\tau+6})^k  \epsilon_0^{k} \leq 1\ ,
\label{eq:eps^2par}
\end{equation}
and holds true provided
\begin{equation}
\epsilon_0 \leq\frac{1}{\Lambda  \alpha^{3\tau+6}}\ .
\label{eq:eps0-1}
\end{equation}

Consider now the restrictions $\delta_k$.  We immediately get
\begin{equation}
\sum_{k\geq1} \delta_k = \sum_{k\geq1} \alpha^{-k}\leq\frac{1}{8}\ ,\quad
\hbox{for }\alpha\geq9\ .
\label{eq:alpha_condition}
\end{equation}
We now prove that~\eqref{eq:sigma*} and~\eqref{eq:rho*} hold true. Starting with
$$
\ln\prod_{k\geq1}(1-4\delta_k) = \sum_{k\geq1} \ln(1-4\delta_k)\ ,
$$
we easily get
$$
0\geq\sum \ln(1-4\delta_k) \geq -8\ln2 \sum_{k\geq1}\delta_k>-\ln2\ ,
$$
from which we have $\rho^*=\rho/2$ and $\sigma^*=\sigma/2$.

Let us now focus on the sequence of the detuning
frequency vectors $\{{\delta\omega}^{(k)}\}_{k\geq0}$, which requires some
additional care. Indeed, Lemma~\ref{lem:iterativo} holds true provided
$\|\delta\omega^{(k)}\cdot p\|_{1-\delta_k} \leq \frac{\epsilon_{k-1}}{2 \delta_k^{\tau+2}}$ and the sequence $\epsilon_k$, by definition, satisfy $\lim_{k\to\infty}{\delta\omega}^{(k)}=0$. The recursive
definition in~\eqref{frm:hat{H}^{(r)}-functions} allows us to compute
${\delta\omega}^{(k)}\cdot{p}$ as
\begin{equation}
{\delta\omega}^{(k)}\cdot{p}=\sum_{j\geq k+1}
\big({\delta\omega}^{(j-1)}\cdot{p}-{\delta\omega}^{(j)}\cdot{p}\big)\ ,
\quad\hbox{for } k \ge 0\ ,
\label{eq:per-delta-omega-come-serie}
\end{equation}
and by using the
inequality~(\ref{frm:stima-detunings}), we get
\begin{equation*}
 \|\delta\omega^{(k)}\cdot p\|_{1-\delta_k} \leq \sum_{j=k+1}^{\infty}  \frac{K_5}{\delta_j^{\tau+2}} \epsilon_{j-1} \ .
\end{equation*}
Thus, the applicability of the Iterative Lemma~\ref{lem:iterativo} is then verified \emph{a posteriori}, if the inequality
\begin{equation}
  \sum_{j\geq k+1} \frac{K_5}{\delta_j^{\tau+2}} \epsilon_{j-1}
  \leq
  \frac{\epsilon_{k-1}}{\delta_k^{\tau+2}}\ 
  \label{eq:aprioridetuning} 
\end{equation}
holds true for every positive integer $k$. We can rewrite this condition as
$$
K_5\sum_{j\geq k+1} (\alpha^{\tau+2} \epsilon_0)^j
=
K_5 \frac{(\alpha^{\tau+2} \epsilon_0)^{k+1}}{1-\alpha^{\tau+2} \epsilon_0}
\leq
(\alpha^{\tau+2} \epsilon_0)^k
$$
from which we get
\begin{equation}
\epsilon_0
\leq
\frac{1}{\alpha^{\tau+2}(K_5+1)}\ .
\label{eq:eps0-2}
\end{equation}

Hence, once the choice of $\alpha$ is made so as to satisfy
\eqref{eq:alpha_condition}, $\alpha\geq9$, one has two additional
smallness conditions on $\epsilon_0$, \eqref{eq:eps0-1}
and~\eqref{eq:eps0-2}, that affects the threshold on the small
parameter $\epsilon^*$.

It remains to prove that the canonical transformation is well defined
on some domain.  To this end, consider the sequence of domains
$\{\Delta_{\rho_k,\sigma_k}\}_{k\geq0}$ with $\rho_k$ and $\sigma_k$
as in~\eqref{eq:rho_k} and~\eqref{eq:sigma_k}.

Then the canonical transformation
$\hat\Cscr^{(k)}:\Delta_{\rho_k,\sigma_k}\to\Delta_{\rho_{k-1},\sigma_{k-1}}$
is analytic. Therefore, by composition, the transformation
$$
\Cscr^{(k)}=\hat\Cscr^{(k)}\circ\cdots\circ\hat\Cscr^{(1)}
$$
is canonical and analytic. Moreover, in view of~\eqref{stima-p} and~\eqref{stima-q} we have
$$
|p^{(k)}-p^{(k-1)}|\leq \sigma \sum_{j=1}^{k}\delta_j^{\tau+3}
\quad\hbox{and}\quad
|q^{(k)}-q^{(k-1)}|\leq \rho \sum_{j=1}^{k}\delta_j^{\tau+3}\ ,
$$
thus, since $\sum_{j\geq1}\delta_j$ is convergent, the sequence
$\{\Cscr^{(k)}\}_{k\geq1}$ converges absolutely to
$$
\Cscr^{(\infty)} : \Delta_{\rho_*,\sigma_*}\to \Delta_{\rho_0,\sigma_0}\ ,
$$
with $\rho_*=\rho_0/2$ and $\sigma_*=\sigma_0/2$.  The absolute
convergence implies the uniform convergence in any compact subset of
$\Delta_{\rho_*,\sigma_*}$, hence $\Cscr^{(\infty)}$ is analytic.
Finally, denoting by $(p^{(\infty)},q^{(\infty)})$ the canonical
coordinates in $\Delta_{\rho_*,\sigma_*}$, and we immediately get
$$
|p_j^{(\infty)}-p_j^{(0)}|\leq \frac{\sigma}{8^{\tau+3}}
\quad\hbox{and}\quad
|q_j^{(\infty)}-q_j^{(0)}|\leq \frac{\rho}{8^{\tau+3}}\ .
$$

Lastly, we now focus on the sequence of detunings. We can bound the norm of
$\delta\omega^{(0)}$ exploiting the recursive definition
$$
\delta\omega^{(0)}\cdot p = \sum_{j\geq 1} \left(\delta\omega^{(j-1)}\cdot p - \delta\omega^{(j)}\cdot p\right)\ .
$$
Indeed, one easily gets
$$
\|\delta\omega^{(0)}\cdot p\|_{\frac{1}{2}}
\leq \sum_{j\geq 1}
\frac{\Lambda}{\delta_j^{\tau+2}}\frac{\epsilon_{j-1} E}{2}
\leq
\frac{E}{2}
\sum_{j\geq 1}
\frac{\Lambda \epsilon_{j-1}}{\delta_{j}^{\tau+2}}
\leq
\frac{E}{2}
\sum_{j\geq 1}
\delta_{j}^{2\tau+4}
\leq\frac{E}{2}\frac{1}{8^{2\tau+4}}
\ .
$$
By the properties of the Lie series transformation, one also has that
the sequence $\{H^{(k)}\}_{k\geq0}$ converges to an analytic function
$H^{(\infty)}$ which by construction is in normal form. This concludes
the proof of Proposition \ref{pro:quantitativa}.

\section*{Acknowledgments}
The authors have been partially supported by the MIUR-PRIN 20178CJA2B
``New Frontiers of Celestial Mechanics: theory and Applications'', by the
MIUR Excellence Department Project awarded to the Department of
Mathematics of the University of Rome ``Tor Vergata'' (CUP
E83C18000100006) and
by the National Group of Mathematical Physics (GNFM-INdAM).

\end{document}